\documentclass[12pt]{article}
\usepackage{hyperref}
\hypersetup{pdfstartview=FitH}
\hypersetup{pdfpagemode=FitWidth} 
\usepackage{times}
\usepackage[latin1]{inputenc}
\newcommand{\n}{\newcommand}
\n{\bb}{\bigskip}
\n{\cl}{\centerline}
\n{\ii}{\hskip5mm\relax}
\begin{document}

{\footnotesize For the last version of this text, type {\em hodgegaillard} on Google. Date of this version: Mon Dec 24 09:01:37 CET 2007.}\vfill

\cl{\Large A Hodge Theorem for Noncompact Manifolds} \vfill

{\bf Theorem}\ \ {\it If $M$ is a riemannian manifold, then the inclusion of the complex of coclosed harmonic forms into the de Rham complex induces a linear isomorphism in cohomology. If $M$ has at most countably many connected components, this linear isomorphism is a Fr\'echet isomorphism.}\bb

[Manifolds are assumed to be $C^\infty$ and Hausdorff.] \bb

Theorem 5 in Section I.9.10 of Bourbaki [2] implying that $M$ is paracompact, we can assume that it is connected, and also that it is non-compact (the result being classical in the compact case). Then the claim follows easily (using the Open Mapping Theorem and the fact that the de Rham cohomology is a Fr\'echet space) from the surjectivity of the laplacian on the de Rham complex. Let us check this surjectivity. In [4, p.~158] de Rham proves (using results of Aronszajn, Krzywicki and Szarski [1]) that a harmonic form which has a zero of infinite order vanishes identically; this implies in particular that the laplacian satisfies Property (A) in Definition~5 of Malgrange [3, p.~333]; it is well known that the laplacian satisfies also Condition (P) --- called {\bf ellipticity} nowadays --- in Definition 6 of [3, p.~338];  in view of Theorem 5 in [3, p.~341] this implies the desired surjectivity. 

\begin{itemize}
\item[{[1]}] Aronszajn N., Krzywicki A., Szarski J., 
A unique continuation theorem for exterior differential forms on Riemannian manifolds, {\it Ark. Mat.}, Volume~4, Number~5 (1962) 417-453. Available at 
\href{http://www.springerlink.com/content/c853518221310827/?p=5beb9fd076f6487fb274a999415eb402&pi=5}%
{springerlink.com}.
\item[{[2]}] Bourbaki, N., {\bf Topologie g\'en\'erale}, Vol. 1, Chapitres 1 \`a 4, Hermann, 1971.
\item[{[3]}] Malgrange B., Existence et approximation des solutions des équations aux dérivées partielles et des équations de convolution, {\it Ann. Inst. Fourier}, Grenoble {\bf 6} (1955-56) 271-354. Available at
\href{http://www.numdam.org/numdam-bin/fitem?id=AIF_1956__6__271_0}%
{numdam.org}.
\item[{[4]}] de Rham G., {\bf Differentiable manifolds}, Springer-Verlag, 1984. 
\end{itemize}
\end{document}